\documentclass[10pt,a4paper,twoside]{article}
\usepackage{amsmath,amsfonts,amsthm,amsopn,color,amssymb,enumitem}
\usepackage{palatino}
\usepackage{graphicx}
\usepackage[colorlinks=true]{hyperref}
\hypersetup{urlcolor=blue, citecolor=red, linkcolor=blue}

\newcommand{\eps}{\varepsilon}

\newcommand{\R}{\mathbb{R}}

\newcommand{\Rn}{{\mathbb{R}^n}}
\newcommand{\RN}{{\mathbb{R}^N}}
\newcommand{\RT}{{\mathbb{R}^3}}

\newcommand{\de}{\partial}

\renewcommand{\le}{\leqslant}
\renewcommand{\ge}{\geqslant}
\renewcommand{\a }{\alpha }

\renewcommand{\b }{\beta }

\renewcommand{\d }{\delta }

\newcommand{\g }{\gamma }

\renewcommand{\l }{\lambda}
\newcommand{\n }{\nabla }
\newcommand{\s }{\sigma }
\renewcommand{\t}{\theta}
\renewcommand{\O}{\Omega}

\newcommand{\G}{\Gamma}

\newcommand{\Itq}{I^T_{q}}

\newcommand{\Iq}{I_{q}}

\renewcommand{\H}{H^1(\RN)}
\newcommand{\Hr}{H^1_r(\RN)}
\newcommand{\HT}{H^1(\RT)}
\newcommand{\HTr}{H^1_r(\RT)}

\newcommand{\E}{\mathcal{E}}

\newcommand{\N}{\mathbb{N}}

\newcommand{\D }{{\mathcal D}^{1,2}(\RT)}
\newcommand{\irn }{\int_{\RN}}
\newcommand{\irt }{\int_{\RT}}

\def\bbm[#1]{\mbox{\boldmath $#1$}}
\newcommand{\beq }{\begin{equation}}
\newcommand{\eeq }{\end{equation}}
\newcommand{\SU}{\sum_{i=1}^k}

\newtheorem{theorem}{Theorem}[section]
\newtheorem{lemma}[theorem]{Lemma}

\newtheorem{proposition}[theorem]{Proposition}
\newtheorem{remark}[theorem]{Remark}

\renewenvironment{proof}{\noindent{\textbf{Proof\quad}}}{$\hfill\square$\vspace{0.2 cm}\\}
\newenvironment{proofmain}{\noindent{\textbf{Proof of Theorem  \ref{main}\quad}}}{$\hfill\square$\vspace{0.2 cm}\\}

\title{{\bf
Multiple critical points\\for a class of nonlinear functionals
\footnote{The authors are supported by M.I.U.R. -
P.R.I.N. ``Metodi variazionali e topologici nello studio di
fenomeni non lineari''}}}

\author{A. Azzollini \thanks{Dipartimento di Matematica ed Informatica, Universit\`a degli
Studi della Basilicata,  Via dell'Ateneo Lucano 10, I-85100
Potenza, Italy, e-mail: {\tt antonio.azzollini@unibas.it}}
 \; \& \;
P. d'Avenia\thanks{Dipartimento di Matematica, Politecnico di
Bari, Via E. Orabona 4, I-70125 Bari, Italy, e-mail: {\tt
p.davenia@poliba.it}}
 \; \& \;
A. Pomponio\thanks{Dipartimento di Matematica, Politecnico di
Bari, Via E. Orabona 4, I-70125 Bari, Italy, e-mail: {\tt
a.pomponio@poliba.it}}}

\date{}

\begin{document}

\maketitle

\begin{abstract}
In this paper we prove a multiplicity result concerning the critical points of a class of
functionals involving local and nonlocal nonlinearities.
We apply our result to the nonlinear Schr\"odinger-Maxwell system in $\RT$
and to the nonlinear elliptic Kirchhoff equation in $\RN$ assuming
on the \emph{local} nonlinearity the general hypotheses introduced by
Berestycki and Lions.
\end{abstract}

\section{Introduction}

In the celebrated papers \cite{BL1,BL2}, Berestycki and Lions proved the existence of a
ground state and a multiplicity result for the equation
    \begin{equation}\label{eq:BL}
        -\Delta u = g(u), \; u:\RN\to\R,
    \end{equation}
for $N\ge 3$, assuming that
\begin{enumerate}[label=(g\arabic*), ref=(g\arabic*)]
  \item \label{g1} $g\in C(\R,\R)$ and odd;
  \item \label{g2} $-\infty <\liminf_{s\to 0^+} g(s)/s\le \limsup_{s\to 0^+} g(s)/s=-m<0$;
  \item \label{g3} $-\infty \le\limsup_{s\to +\infty} g(s)/s^{2^*-1}\le 0$, with $2^*=2N/(N-2)$;
  \item \label{g4} there exists $\zeta>0$ such that $G(\zeta):=\int_0^\zeta g(s)\,d  s>0$.
\end{enumerate}
Modifying, if necessary, in a suitable way the nonlinearity $g$
(without losing the generality of the problem), it can be proved
that equation \eqref{eq:BL} possesses a variational structure,
namely its solutions can be found as critical points of the
functional
\[
I(u) = \frac 1 2\irn|\n u|^2 - \irn G(u).
\]
Solutions of several nonlinear elliptic equations involving local and nonlocal nonlinearities can be found looking for critical points of a suitable perturbation of $I$, namely
    \begin{equation}\label{eq:funq}
        \Iq(u)=\frac{1}{2} \irn |\nabla u|^2+q R(u)-\irn G(u),\;u\in\H,
    \end{equation}
where $q>0$ is a small parameter and $R:\H\to\R$. In order to define the functional $I_q$ we need to replace \ref{g3} with the stronger assumption
\begin{enumerate}[label=(g\arabic*)', ref=(g\arabic*)']\setcounter{enumi}{2}
  \item \label{g3'} $\lim_{s\to + \infty} g(s)/|s|^{2^*-1}= 0$.
\end{enumerate}

In this paper we are interested in providing a multiplicity result in critical point theory for $I_q$.
To this end we suppose that $R=\sum_{i=1}^k R_i$ and, for each $i=1,\ldots ,k$ the functional $R_i$ satisfies:
\begin{enumerate}[label=(R\arabic*), ref=(R\arabic*)]
    \item \label{r1} $R_i$ is $C^1(\H,\R)$, nonnegative and even;
    \item \label{rd'} there exists $\d_i>0$ such that $R_i'(u)[u]\le C\|u\|^{\d_i}$, for any $u \in \H$;
    \item \label{r4} if $\{u_j\}_j$ is a sequence in $\H$ weakly convergent to $u\in\H$,
    then
    \[
    \limsup_j R_i'(u_j)[u-u_j]\le 0;
    \]
    \item \label{r5} there exist $\a_i,\b_i\ge 0$ such that if $u\in\H$, $t>0$ and $u_t=u(\cdot /t)$, then
    \[
    R_i(u_t)=t^{\a_i} R_i(t^{\b_i} u);
    \]
    \item \label{r6} $R_i$ is invariant under the action of $N$-dimensional orthogonal group,
    i.e. $R_i(u(g\cdot))=R_i(u(\cdot))$ for every $g\in O(N)$.
\end{enumerate}


The effect deriving from the presence of the perturbation
$qR$ is to modify the structure of the functional $I$ both as regards
the geometrical properties, and as regards compactness properties.
In particular two remarkable difficulties arise: the first is
related with the problem of applying classical min-max arguments to
find Palais-Smale sequences at suitable levels, the second is
concerned with the compactness of these sequences. If, on one hand,
just assuming the positiveness of the functional $R$ we overcome
the difficulty of finding suitable min-max levels, on the other, the
problem of boundedness of Palais-Smale sequences is not nearly
trivial. This is a consequence of the fact that no
Ambrosetti-Rabinowitz hypothesis like
    \begin{equation*}
        0< \nu G(t)\le
        tg(t), \;\hbox{for } \nu >2,
    \end{equation*}
is assumed on $g$. The monotonicity trick based on an idea of Struwe
\cite{Sw} and formalized by Jeanjean \cite{J} has turned out to be a
powerful method to overcome this difficulty. By means of the
monotonicity trick and a truncation argument based on an idea of Berti and Bolle \cite{BB} and of 
Jeanjean and Le Coz \cite{JC} (see also \cite{K2}), in \cite{ADP}
we have proved an existence result for a functional which
is included in the class we are treating. The same arguments have
been used also in \cite{A} to prove a similar existence result also
for another functional of the type described in \eqref{eq:funq}. In
both the results it is required that the parameter $q$ is
sufficiently small. The well known fact proved in \cite{BL2} and
more recently in \cite{HIT} that $I$ possesses infinitely many
critical points has led us to wonder if, at least for small $q,$ a
multiplicity result on the number of critical points keeps holding
for $I_q$. In this direction a fundamental contribution comes from
the recent paper \cite{HIT}, where, developing some ideas of \cite{J0}, a new method to find multiple
solutions to equations involving general local nonlinearities has
been introduced. Here we will get our multiplicity result by a
suitable combination of the new method described in \cite{HIT} with
the
truncation argument of \cite{JC}.\\

Our main result is the following.

\begin{theorem}\label{main}
Let us suppose \ref{g1}, \ref{g2}, \ref{g3'}, \ref{g4} and \ref{r1}--\ref{r6}. Then for
any $h\in \N,h\ge 1$, there exists $q(h)>0$ such that for any $0<q<q(h)$ the functional $I_q$
admits at least $h$ couples of critical points in $\H$ with radial symmetry.
\end{theorem}

Some nonlinear mathematical physics problems can be solved looking for critical points of
functionals strictly related with $I_q$. Among them, we recall, for instance,
the electrostatic Schr\"odinger-Maxwell equations.
This system constitutes a model to describe the interaction between a nonrelativistic
charged particle with the electromagnetic field (see for example \cite{AR,A0,ADP,AP08,BF,C2,DM1,DA,JZ,K1,K2,Ru,WZ,ZZ}). In the electrostatic case the system becomes
    \begin{equation}\label{eq:sm}
        \left\{
            \begin{array}{ll}
                -\Delta u+q\phi u=g(u)&\hbox{in }\RT,
                \\
                -\Delta \phi=q u^2&\hbox{in }\RT.
            \end{array}
        \right.
    \end{equation}
Finding solutions to the previous system is equivalent to look for
critical points of the functional
\begin{equation*}
I_q(u)= \frac 12 \irt |\n u|^2 + \frac q4 \irt \left(\frac 1 {|x|}\ast u^2\right)u^2
-\irt G(u).
\end{equation*}
In \cite{AR}, the authors study \eqref{eq:sm} with $g(u)=-u+|u|^{p-1}u$ and $1<p<5$ and use an abstract tool,
based on the monotonicity trick, to prove a multiplicity result.

As a consequence of Theorem \ref{main} we prove

\begin{theorem}\label{main2}
Let us suppose \ref{g1}, \ref{g2}, \ref{g3}, \ref{g4}. Then for
any $h\in \N,h\ge 1$, there exists $q(h)>0$ such that for any $0<q<q(h)$ system \eqref{eq:sm}
admits at least $h$ couples of solutions in $\HT\times\D$ with radial symmetry.
\end{theorem}

Another variational problem related with our abstract result is the following.
Let us consider the multidimensional Kirchhoff equation
\[
\frac{\partial^2 u}{\partial t^2} - \left(p+q\int_\O |\n u|^2\right) \Delta u=0
\qquad \hbox{ in }\O,
\]
where $\O\subset\R^N$, $p>0$ and $u$ satisfies some initial or
boundary conditions. It arises from the following Kirchhoff'
nonlinear generalization (see \cite{K}) of the well known d'Alembert
equation
    \begin{equation*}
        \rho\frac{\partial^2 u}{\partial t^2} - \left(\frac{P_0}{h}+\frac {E}{2L}
        \int_0^L \left|\frac{\partial u}{\partial x}\right|^2\,dx
        \right) \frac{\partial^2 u}{\partial x^2}=0,
    \end{equation*}
and it describes a vibrating string, taking into account the changes
in length of the string during the vibration. Here, $L$ is the
length of the string, $h$ is the area of the cross section, $E$ is
the Young modulus of the material, $\rho$ is the mass density and
$P_0$ is the initial tension.\\
If we look for static solutions, the equation we have to solve is
    \begin{equation*}
        - \left(p+q\int_\O |\n u|^2\right) \Delta u=0.
    \end{equation*}
In the same spirit of \cite{ACM,A} we consider the semilinear perturbation
    \begin{equation}\label{K}
        - \left(p+q\int_\O |\n u|^2\right) \Delta u=g(u),\quad \hbox{in }\O\subset \RN.
    \end{equation}
Recently this equation has been extensively treated by many authors
in bounded domains, assuming Dirichlet conditions on the boundary
(see for example \cite{ACM,HZ,M,MZ,PZ,Ri,PZ2}).\\
Here we are
interested in showing an application of our abstract result to the
equation \eqref{K} in all the space $\RN$, $N\ge 3$. The solutions are the
critical points of the functional
    \begin{equation*}
        I_q(u)= \frac p 2 \irn|\n u|^2 + \frac{q}{4} \left(  \irn|\n u|^2\right)^2-\irn G(u).
    \end{equation*}
We prove the following result.
\begin{theorem}\label{main3}
Let us suppose \ref{g1}, \ref{g2}, \ref{g3}, \ref{g4}. Then for
any $h\in \N,h\ge 1$, there exists $q(h)>0$ such that for any $0<q<q(h)$
equation \eqref{K} admits at least $h$ couples of solutions in $\H$ with radial symmetry.
\end{theorem}

The paper is organized as follows: in Section \ref{abs} we prove
Theorem \ref{main}; in Section \ref{app} we show as it can be
applied to the nonlinear Schr\"odinger-Maxwell system and the
nonlinear elliptic Kirchhoff equation in order to prove Theorems \ref{main2} and \ref{main3}.

\begin{center}
{\bf NOTATION}
\end{center}
We will use the following notations:
\begin{itemize}
\item for any $1\le s\le +\infty$, we denote by $\|\cdot\|_s$ the usual norm of the Lebesgue space $L^s(\RN)$;
\item $H^1(\RN)$ is the usual Sobolev space endowed with the norm
\[
\|u\|^2:=\irn |\n u|^2+ u^2;
\]
\item ${\cal D}^{1,2}(\RN)$ is completion of $C_0^\infty(\RN)$ (the compactly
supported functions in $C^\infty(\RN)$) with respect to the norm
\[
\|u\|_{{\cal D}^{1,2}(\RN)}^2:=\irn |\n u|^2;
\]
\item $C,C',C_i$ are various positive constants which may also vary from line to line.
\end{itemize}

\section{The abstract result}\label{abs}

We set for any $s\ge 0$,
\begin{align*}
g_1(s) &:=
(g(s)+ms)^+,
\\
g_2(s) &:=g_1(s)-g(s),
\end{align*}
and we extend them as odd functions.
Since
    \begin{align}
        \lim_{s\to 0}\frac{g_1(s)}{s} &= 0,\nonumber
        \\
        \lim_{s\to\pm\infty}\frac{g_1(s)}{|s|^{2^*-1}}&=0,\label{eq:lim2}
    \end{align}
and
    \begin{equation}
        g_2(s) \ge ms,\quad  \forall s\ge 0,\label{eq:g2}
    \end{equation}
by some computations, we have that for any $\eps>0$ there exists
$C_\eps>0$ such that
    \begin{equation}
        g_1(s) \le C_\eps |s|^{2^*-1}+\eps g_2(s),\quad  \forall
        s\ge0\label{eq:g1g2}.
    \end{equation}
If we set
    \begin{equation*}
        G_i(t):=\int^t_0g_i(s)\,ds,\quad i=1,2,
    \end{equation*}
then, by \eqref{eq:g2} and \eqref{eq:g1g2}, we have
    \begin{equation}
         G_2(s) \ge \frac m 2 s^2,\quad  \forall s\in\R\label{eq:G2}
    \end{equation}
and for any $\eps>0$ there exists $C_\eps>0$ such that
    \begin{equation}
        G_1(s) \le C_\eps |s|^{2^*}+\eps G_2(s),\quad  \forall
        s\in\R\label{eq:G1G2}.
    \end{equation}

Since, for any $u\in \H$, $R_i(u)-R_i(0)=\int_0^1 \frac{d}{d t}R_i(t u) d t$, by \ref{rd'} we have that
\begin{equation}\label{rd}
R_i(u)\le C_1 +C_2\|u\|^{\d_i}.
\end{equation}

The hypothesis \ref{r6} assures that all functionals that we
will consider in this paper are invariant under rotations. Then
\[
\Hr=\{u\in \H\mid u \hbox{ radial }\}
\]
is a natural constraint to look for critical points,
namely critical points of the functional restricted to $\Hr$
are \emph{true} critical points in $\H$.
Therefore, from now on, we will directly define our functionals in $\Hr$.

As in \cite{JC}, we consider a cut-off function $\chi\in C^\infty(\R_+,\R)$ such that
\begin{equation*}
    \left\{
\begin{array}{ll}
    \chi(s)=1,&\hbox{for }s\in[0,1],\\
    0\le \chi(s)\le 1,&\hbox{for }s\in]1,2[,\\
    \chi(s)=0,&\hbox{for }s\in[2,+\infty[,\\
    \|\chi '\|_\infty \le 2,&
\end{array}
    \right.
\end{equation*}
and we introduce the following truncated functional $\Itq :\Hr\to \R$
\begin{equation*}
\Itq(u)=\frac 12 \irt |\n u|^2 + q k_T (u) R(u)  - \irt G(u),
\end{equation*}
where
\[
k_T(u)=\chi\left(\frac{\| u\|^2}{T^2} \right).
\]
Of course, any critical point $u$ of $\Itq$ with $\|u\|\le T$ is a critical point of $I_q$.

The $C^1-$functional $\Itq$ has the symmetric mountain pass geometry:
\begin{lemma}\label{le:gon}
There exist $r_0>0$ and $\rho_0>0$ such that
\begin{align}
\Itq(u)\ge 0, &\quad\hbox{for } \|u\|\le r_0,    \label{eq:gon1}\\
\Itq(u)\ge\rho_0, &\quad\hbox{for } \|u\|=r_0.   \label{eq:gon2}
\end{align}
Moreover, for any $n\in\N, n\ge 1$, there exists an odd continuous map
\begin{equation*}
\g_{n}:S^{n-1}=\{\s=(\s_1,\cdots,\s_n)\in \Rn\mid |\s|=1\}\to \Hr,
\end{equation*}
such that
\[
\Itq(\g_{n}(\s))<0,\qquad \hbox{for all }\s\in S^{n-1}.
\]
\end{lemma}

\begin{proof}
By \eqref{eq:G2}, \eqref{eq:G1G2} and the positivity of the map $R$,
\[
\Itq (u)\ge C_1 \|u\|^2- C_2 \|u\|^{2^*}
\]
from which we obtain \eqref{eq:gon1} and \eqref{eq:gon2}.\\
Moreover, arguing as in \cite[Theorem 10]{BL2}, for every $n\ge 1$ we can consider
an odd continuous map $\pi_n : S^{n-1}\to\Hr$ such that
\[
0\notin \pi_n (S^{n-1}),\qquad \irn G(\pi_n(\s)) \ge 1 \hbox{ for all } \s\in S^{n-1}.
\]
Then, for $t$ sufficiently large, we take
\[
\g_n(\s)=\pi_n(\s)(\cdot /t)
\]
and we obtain
\begin{align*}
\Itq(\g_n(\s))
&=\frac{t^{N-2}}{2} \irn |\nabla \pi_n(\s)|^2
+q\chi \left( \frac{t^{N-2}\|\n \pi_n(\s)\|_2^2+t^{N}\|\pi_n(\s)\|_2^2}{T^2} \right)
R(\g_n(\s))\\
&\quad-t^{N}\irn G(\pi_n(\s))\\
&\le \frac{t^{N-2}}{2} \irn |\nabla \pi_n(\s)|^2 - t^{N} <0.
\end{align*}
\end{proof}

Let us define
\[
b_n=b_n(q,T)=\inf_{\gamma\in\Gamma_n} \max_{\sigma\in D_n} \Itq(\gamma(\sigma))
\]
where $D_n=\{\sigma=(\s_1,\cdots,\s_n)\in\Rn   \;\vline\;|\sigma|\le 1\}$,
\[
\Gamma_n=\left\{\gamma\in C(D_n, \Hr)\;\vline\;
\begin{array}{ll}
    \gamma(-\sigma)=-\gamma(\sigma)    &\hbox{for all }\sigma\in D_n\\
    \gamma(\sigma)=\gamma_{n}(\sigma) &\hbox{for all }\sigma\in \partial D_n
\end{array}
\right\}
\]
and $\gamma_{n}:\partial D_n \to \Hr$ is given in Lemma \ref{le:gon}.

Analogously to \cite{HIT}, we set
\begin{align*}
\tilde\Iq (\t,u)&=\Iq(u(e^{-\t}\cdot)),
\\
\tilde\Itq (\t,u)&=\Itq(u(e^{-\t}\cdot)),
\\
\tilde I'_q (\t,u)&=\frac{\de }{\de u}\tilde \Iq(\t,u),
\\
(\tilde \Itq)'(\t,u)&=\frac{\de }{\de u}\tilde \Itq(\t,u),
\\
\tilde b_n=\tilde b_n (q,T)&=\inf_{\tilde\gamma\in\tilde\Gamma_n} \max_{\sigma\in D_n} \tilde \Itq(\tilde \gamma(\sigma)),
\end{align*}
where
\[
\tilde\Gamma_n=\left\{\tilde\gamma\in C(D_n,\R\times \Hr)\;\vline\;
\begin{array}{ll}
  \tilde\gamma(\sigma)=(\t(\sigma),\eta(\sigma)) \hbox{ satisfies}   & \\
    (\t(-\sigma),\eta(-\sigma))=(\t(\sigma),-\eta(\sigma))    &\hbox{for all }\sigma\in D_n\\
    (\t(\sigma),\eta(\sigma))=(0,\gamma_{n}(\sigma))         &\hbox{for all }\sigma\in \partial D_n
\end{array}
\right\}.
\]
By \ref{r5} we have
\begin{align*}
\tilde{\Iq}(\t,u)&=\frac{e^{(N-2)\t}}{2} \irn |\nabla u|^2+q \SU e^{\a_i \t}R_i(e^{\b_i\t}u)-e^{N\t}\irn G(u),\\
\tilde\Itq (\t,u)&
=\frac{e^{(N-2)\t}}{2} \irn |\nabla u|^2
+q\chi \left( \frac{e^{(N-2)\t}\|\n u\|_2^2+e^{N\t}\|u\|_2^2}{T^2} \right)
\SU e^{\a_i\t}R_i(e^{\b_i\t}u)
-e^{N\t}\irn G(u).
\end{align*}
Arguing as in \cite{HIT}, the following lemmas hold.
\begin{lemma}\label{le:bn}
We have
\begin{enumerate}
    \item there exists $\bar b>0$ such that $b_n\ge \bar b$, for any $n\ge 1$;
    \item $b_n\to +\infty$;
    \item $b_n=\tilde b_n$, for any $n\ge 1$.
\end{enumerate}
\end{lemma}

\begin{lemma}\label{pr:3.2}
For any $n\ge 1$, there exists a sequence $\{(\t_j,u_j)\}_j\subset \R \times \Hr$ such that
\begin{enumerate}[label=(\roman*), ref=(\roman*)]
    \item $\t_j\to 0$;
    \item $\tilde \Itq(\t_j,u_j)\to b_n$;
    \item $(\tilde \Itq)'(\t_j,u_j)\to 0$ strongly in $(\Hr)^{-1}$; \label{eq:iii3.2}
    \item $\frac{\de }{\de \t}\tilde \Itq(\t_j,u_j)\to 0$.
\end{enumerate}
\end{lemma}

Now we prove that for a suitable choice of $T$ and $q$, the sequence $\{u_j\}_j$ obtained in Lemma \ref{pr:3.2}
actually is a bounded Palais-Smale sequence fot $I_q$.

\begin{proposition}\label{th:T}
Let $n\ge 1$ and $T_n>0$ sufficiently large. There exists $q_n$ which depends on $T_n$, such that for any $0<q<q_n$,
if $\{(\t_j,u_j)\}_j\subset \R \times \Hr$ is the sequence given in Lemma \ref{pr:3.2}, then, up to a subsequence, $\|u_j\|\le T_n$, for any $j\ge 1$.
\end{proposition}

\begin{proof}
By Lemmas \ref{le:bn} and \ref{pr:3.2}, we infer that
\[
N \tilde \Itq(\t_j,u_j)- \frac{\de}{\de \t}\tilde \Itq(\t_j,u_j)=N b_n+o_j(1),
\]
and so
\begin{align}
e^{(N-2)\t_j}\irn|\n u_j|^2
&=q\chi\left( \frac{\|u_j(e^{-{\t_j}}\cdot)\|^2}{T^2} \right)
\SU(\a_i-N) R_i(u_j(e^{-{\t_j}}\cdot))\nonumber
\\
&\quad+q \chi\left( \frac{\|u_j(e^{-{\t_j}}\cdot)\|^2}{T^2} \right)
\SU e^{\a_i \t_j}R'_i(e^{\b_i\t_j}u_j)[\b_i e^{\b_i \t_j}u_j]\nonumber
\\
&\quad+q  \chi'\left( \frac{\|u_j(e^{-{\t_j}}\cdot)\|^2}{T^2} \right)
\frac{(N-2)e^{(N-2)\t_j}\|\n u_j\|_2^2+Ne^{N\t_j}\|u_j\|_2^2}{T^2}
R(u_j(e^{-{\t_j}}\cdot))\nonumber
\\
&\quad+N b_n+o_j(1).\label{eq:ineq}
\end{align}
We are going to estimate the right part of the previous identity.
By the min-max definition of $b_n$, if $\g \in \G_n$, we have
\begin{align*}
b_n & \le
\max_{\s\in D_n}\Itq(\g(\s))
\\
& \le \max_{\s \in D_n}\left\{\frac 1 2 \irn |\n \g(\s)|^2
-\irn G(\g(\s))\right\}
 +\max_{\s\in D_n}
\left\{q k_T(\g(\s)) R(\g(\s))\right\}
\\
& = A_1 + A_2(T)
\end{align*}
Now, if $\|\g(\s)\|^2\ge 2 T^2$ then $A_2(T)=0,$ otherwise, by \eqref{rd}, we have
\begin{align*}
A_2(T)\le  q (C_1+C_2\|\g(\s)\|^\d)\le  q (C_1+C'_2 T^\d),
\end{align*}
for a suitable $\d>0$. Moreover we have that
\begin{align}
&q\chi\left( \frac{\|u_j(e^{-{\t_j}}\cdot)\|^2}{T^2} \right)
\SU(\a_i-N) R_i(u_j(e^{-{\t_j}}\cdot))
\le q (C_1+C_2 T^\d);\nonumber
\\
&q \chi\left(\frac{\|u_j(e^{-{\t_j}}\cdot)\|^2}{T^2} \right)
\SU e^{\a_i \t_j}R'_i(e^{\b_i\t_j}u_j)[\b_i e^{\b_i \t_j}u_j]\le C q T^{\d};
\label{eq:st2}
\\
&q  \chi'\left( \frac{\|u_j(e^{-{\t_j}}\cdot)\|^2}{T^2} \right)
\frac{(N-2)e^{(N-2)\t_j}\|\n u_j\|_2^2+Ne^{N\t_j}\|u_j\|_2^2}{T^2}
R(u_j(e^{-{\t_j}}\cdot))\le q (C_1+C_2 T^\d).
\label{eq:st3}
\end{align}
Then, from \eqref{eq:ineq} we deduce that
\begin{equation}\label{eq:first}
\irn |\n u_j|^2 \le C' + q (C_1+C_2 T^\d).
\end{equation}
On the other hand, since $\frac{\de}{\de \t}\tilde \Itq(\t_j,u_j)=o_j(1)$, by \eqref{eq:G1G2} we
have that
\begin{align}
&\frac{(N-2)e^{(N-2)\t_j}}{2}\irn|\n u_j|^2
+q\chi\left( \frac{\|u_j(e^{-{\t_j}}\cdot)\|^2}{T^2} \right)
\SU \a_i R_i(u_j(e^{-{\t_j}}\cdot))         \nonumber
\\
&\quad+q \chi\left( \frac{\|u_j(e^{-{\t_j}}\cdot)\|^2}{T^2} \right)
\SU e^{\a_i \t_j}R'_i(e^{\b_i\t_j}u_j)[\b_i e^{\b_i \t_j}u_j]       \nonumber
\\
&\quad+q  \chi'\left( \frac{\|u_j(e^{-{\t_j}}\cdot)\|^2}{T^2} \right)
\frac{(N-2)e^{(N-2)\t_j}\|\n u_j\|_2^2+Ne^{N\t_j}\|u_j\|_2^2}{T^2}  R(u_j(e^{-{\t_j}}\cdot))      \nonumber
\\
&\quad+N e^{N\t_j}\irn G_2(u_j)
=N e^{N\t_j}\irn G_1(u_j)+o_j(1)  \nonumber
\\
&\le N e^{N\t_j} \left(C_\eps\irn|u_j|^{2^*}
+\eps\irn G_2(u_j)\right)
+o_j(1).\label{eq:neh}
\end{align}
Now, by \eqref{eq:G2}, \eqref{eq:st2}, \eqref{eq:st3}, \eqref{eq:first} and \eqref{eq:neh}, we obtain
\begin{align}\label{eq:second}
&\frac{N e^{N\t_j}m  (1-\eps)}{2}\irn u^2_j  \le (1-\eps)N e^{N\t_j}\irn G_2(u_j)
 \nonumber
\\
&\le N e^{N\t_j}C_\eps \irn |u_j|^{2^*}-q \chi\left( \frac{\|u_j(e^{-{\t_j}}\cdot)\|^2}{T^2} \right)
\SU e^{\a_i \t_j}R'_i(e^{\b_i\t_j}u_j)[\b_i e^{\b_i \t_j}u_j]       \nonumber
\\
&\quad- q  \chi'\left( \frac{\|u_j(e^{-{\t_j}}\cdot)\|^2}{T^2} \right)
\frac{(N-2)e^{(N-2)\t_j}\|\n u_j\|_2^2+Ne^{N\t_j}\|u_j\|_2^2}{T^2}
R(u_j(e^{-{\t_j}}\cdot)) +o_j(1)     \nonumber
\\
&\le C \left(\irn |\n u_j|^2\right)^{2^*/2} + q (C_1+C_2 T^\d) +o_j(1) \nonumber
\\
&\le C(C' + q (C_1+C_2 T^\d))^{2^*/2} + q (C_1+C_2 T^\d)+o_j(1).
\end{align}
We suppose by contradiction that there exists no subsequence of
$\{u_j\}_j$ which is uniformly bounded by $T$ in the $H^1-$norm. As a
consequence, for a certain $j_0$ it should result that
\begin{equation}\label{eq:contr}
\|u_j\|> T, \quad\forall j\ge  j_0.
\end{equation}
Without any loss of generality, we are supposing that
\eqref{eq:contr} is true for any $u_j.$
\\
Therefore, by \eqref{eq:first} and \eqref{eq:second}, we conclude
that
\begin{align*}
T^2 < \|u_j\|^2  \le C_3 +C_4 q T^{\frac{2^*}{2}\d }
\end{align*}
which is not true for $T$ large and $q$ small enough: indeed we can find $T_0>0$ such that $T_0^2>C_3+1$ and $q_0=q_0(T_0)$ such that $C_4 q T^{\frac{2^*}{2}\d }<1$, for any $q<q_0$, and we find a contradiction.
\end{proof}

In our arguments, the following variant of the Strauss'
compactness result \cite{Str}
(see also \cite[Theorem A.1]{BL1})
 will be a fundamental tool.
\begin{proposition}\label{le:str}
Let $P$ and $Q:\R\to\R$ be two continuous functions satisfying
\begin{equation*}
\lim_{s\to\infty}\frac{P(s)}{Q(s)}=0,
\end{equation*}
$\{v_j\}_j,$ $v$ and $w$ be measurable functions from $\RN$ to $\R$,
with $w$ bounded, such that
\begin{align*}
&\sup_j\irn | Q(v_j(x))w|\,dx <+\infty,
\\ 
&P(v_j(x))\to v(x) \:\hbox{a.e. in }\RN. 
\end{align*}
Then $\|(P(v_j)-v)w\|_{L^1(B)}\to 0$, for any bounded Borel set
$B$.

Moreover, if we have also
\begin{align*}
\lim_{s\to 0}\frac{P(s)}{Q(s)} &=0,\\ 
\lim_{|x|\to\infty}\sup_j |v_j(x)| &= 0, 
\end{align*}
then $\|(P(v_j)-v)w\|_{L^1(\RN)}\to 0.$
\end{proposition}

In analogy with the well-known compactness result in \cite{BL2},
we state the following result.

\begin{lemma}\label{le:PS}
Let $n\ge 1$, $T_n,q_n>0$ as in Proposition \ref{th:T} and $\{(\t_j,u_j)\}_j\subset \R \times \Hr$ be the sequence given in Lemma \ref{pr:3.2}. Then $\{u_j\}_j$  admits a subsequence which converges in $\Hr$ to a nontrivial critical point of $I_q$ at level $b_n$.
\end{lemma}

\begin{proof}
Since $\{u_j\}_j$ is bounded, up to a subsequence, we can suppose that there exists $u\in\Hr$
such that
\begin{align}
u_j\rightharpoonup u\;&\hbox{ weakly in }\Hr,\nonumber
\\
u_j\to u\;&\hbox{ in }L^p(\RN),\; 2<p<2^*,\nonumber
\\
u_j\to u\;&\hbox{ a.e. in }\RN.\label{eq:aeconv}
\end{align}
By weak lower semicontinuity we have
\begin{align}
\irn |\n u|^2 \le &\liminf_j \irn |\n u_j|^2. \label{eq:semi1}
\end{align}
Since $\|u_j\|\le T_n$ we have
\begin{align*}
\tilde I'_q(\t_j,u_j)[v]&=(\tilde\Itq)'(\t_j,u_j)[v]\\
&=
e^{(N-2)\t_j}\irn \n u_j \cdot \n v
+q\SU e^{(\a_i+\b_i) \t_j}R'_i(e^{\b_i \t_j} u_j)[v]\\
&\quad
+e^{N\t_j}\irn g_2(u_j)v
-e^{N\t_j}\irn g_1(u_j)v
\end{align*}
for every $v\in\H$.\\
Then, by \ref{eq:iii3.2} of Lemma \ref{pr:3.2}
\begin{align}
&\tilde I'_q(\t_j,u_j)[u]-\tilde I'_q(\t_j,u_j)[u_j]\nonumber
\\
&=e^{(N-2)\t_j}\irn \n u_j\cdot(\n u-\n u_j)
+q\SU e^{(\a_i+\b_i) \t_j}R'_i(e^{\b_i \t_j} u_j)[u-u_j] \nonumber
\\
&\quad+e^{N\t_j}\irn g_2(u_j)(u-u_j)
-e^{N\t_j}\irn g_1(u_j)(u-u_j)
=o_j(1).\label{eq:euj}
\end{align}
If we apply Proposition \ref{le:str} for $P(s)=g_i(s)$, $i=1,2,$
$Q(s)= |s|^{2^*-1},$ $(v_j)_j=(u_j)_j,$ $v=g_i(u),$ $i=1,2$ and
$w$ a generic $C^\infty_0(\RN)$-function, by \ref{g3'}, \eqref{eq:lim2} and
\eqref{eq:aeconv} we deduce that
\[
\irn g_i(u_j)w\to\irn g_i(u)w\quad i=1,2,
\]
and so
\begin{equation}\label{eq:convforteg}
\irn g_i(u_j)u\to\irn g_i(u)u\quad i=1,2.
\end{equation}
Moreover, applying Proposition \ref{le:str} for $P(s)=g_1(s)s,$ $Q(s)= s^2+|s|^{2^*},$ $(v_j)_j=(u_j)_j,$ $v=g_1(u)u,$ and $w=1,$ by
\ref{g3'}, \eqref{eq:lim2} and \eqref{eq:aeconv},
we deduce that
    \begin{align}
        \irn g_1(u_j)u_j \to \irn g_1(u)u.\label{eq:convg}
    \end{align}
Moreover, by \eqref{eq:aeconv} and Fatou's lemma
    \begin{align}
        \irn g_2(u)u\le &\liminf_j \irn g_2(u_j)u_j.\label{eq:convg2}
    \end{align}
By \eqref{eq:euj}, \eqref{eq:convforteg}, \eqref{eq:convg} \eqref{eq:convg2} and \ref{r4}, we have
\begin{align}
\limsup_j \irn  |\n u_j|^2
&=\limsup_j e^{(N-2)\t_j}\irn  |\n u_j|^2 \nonumber
\\
&=\limsup_j \left[e^{(N-2)\t_j}\irn \n u_j\cdot\n u
+q\SU e^{(\a_i+\b_i) \t_j}R'_i(e^{\b_i \t_j} u_j)[u-u_j]\right. \nonumber
\\
&\quad \left.+e^{N\t_j}\irn g_2(u_j)(u-u_j)
-e^{N\t_j}\irn g_1(u_j)(u-u_j)\right] \nonumber
\\
&\le \irn |\n u|^2. \label{eq:limsup}
\end{align}
By \eqref{eq:semi1} and \eqref{eq:limsup}, we get
\begin{align}
\lim_j \irn |\n u_j|^2= &\irn |\n u|^2 \label{eq:sf1},
\end{align}
hence, by \eqref{eq:euj},
\begin{equation}    \label{eq:sfg2}
\lim_j \irn g_2(u_j)u_j= \irn g_2(u)u.
\end{equation}
Since $g_2(s)s=ms^2 +h(s)$, with $h$ a positive and continuous function, by Fatou's Lemma we have
\begin{align*}
\irn h(u) \le &\liminf_j \irn h(u_j),
\\
\irn u^2 \le &\liminf_j \irn u_j^2.
\end{align*}
These last two inequalities and \eqref{eq:sfg2} imply that, up to
a subsequence,
\[
\lim_j \irn u_j^2 = \irn  u^2,
\]
which, together with \eqref{eq:sf1}, shows that $u_j \to u$ strongly in $\Hr$. Therefore, since $b_n>0$, $u$ is a nontrivial critical point of $I_q$ at level $b_n$.
\end{proof}

\begin{proofmain}
Let $h\ge 1$. Since $b_n\to +\infty$, up to a subsequence, we can consider $b_1<b_2<\cdots<b_h$. By Lemma \ref{le:PS} we conclude, defining $q(h)=q_h>0$.
\end{proofmain}

\section{Some applications} \label{app}

\subsection{The nonlinear Schr\"odinger-Maxwell system}

Let us consider the Schr\"odinger-Maxwell system:
\begin{equation}    \label{SM}\tag{${\cal SM}$}
\left\{
\begin{array}{ll}
-\Delta u+q\phi u=g(u)&\hbox{in }\RT,
\\
-\Delta \phi=q u^2&\hbox{in }\RT,
\end{array}
\right.
\end{equation}
where $q>0$ and $g$ satisfies \ref{g1}-\ref{g4}. Arguing as in \cite{ADP,BL1}, without loss of generality, we can suppose that $g$ satisfies \ref{g3'}.

The solutions $(u,\phi)\in \HT \times \D$ of
\eqref{SM} are the critical points of the action functional $\mathcal{E}_q
\colon \HT \times \D \to \R$, defined as
\[
\mathcal{E}_q(u,\phi):=\frac 12 \irt |\n u|^2
-\frac 14 \irt |\n \phi|^2
+\frac q2 \irt \phi u^2
-\irt G(u).
\]
The action functional $\E_q$ exhibits a strong indefiniteness, namely it is
unbounded both from below and from above on infinite dimensional
subspaces. This indefiniteness can be removed using the reduction
method described in \cite{BF}, by which we are led to study a
one variable functional that does not present such a strongly
indefinite nature. Indeed, for every $u\in L^\frac{12}5(\RT)$, there exists a unique $\phi_u\in
\D$ solution of
\[
-\Delta \phi=q u^2,\qquad \hbox{in }\RT.
\]
Moreover it can be proved that $(u,\phi)\in H^1(\RT)\times \D$ is a solution
of \eqref{SM} (critical point of functional $\mathcal{E}_q$) if and
only if $u\in\HT$ is a critical point of the functional $I_q\colon
\HT\to \R$ defined as
\begin{equation*}
I_q(u)= \frac 12 \irt |\n u|^2 + \frac q4 \irt \phi_u u^2
-\irt G(u),
\end{equation*}
and $\phi=\phi_u$.

According to our notations, in this case
$R(u)=\frac 14 \irt \phi_u u^2$. In order to check that $R$
satisfies \ref{r1}-\ref{r6}, we need some preliminary results on $\phi_u$ (see for example \cite{DM1}).
\begin{lemma}\label{le:prop}
The map $u\in L^{\frac{12}{5}}(\RT)\mapsto \phi_u\in \D $ is $C^1$. Moreover, for every $u\in \HT$, we have
\begin{enumerate}[label=\roman*),ref=\roman*)]
\item $\|\phi_u\|^2_{\D}=q\irt\phi_u u^2$;
\item $\phi_u\ge 0$;
\item $\phi_{-u}=\phi_u$;
\item for any $t>0$: $\phi_{u_t}(x)=t^2\phi_u(x/t)$, where $u_t(x)=u(x/t)$;
\item there exist $C,C'>0$ independent of $u\in\HT$ such that
$$\|\phi_u\|_{\D}\le C q \|u\|^2_\frac{12}{5},$$
and
\begin{equation}\label{eq:phiq}
\irt\phi_u u^2\le C'q \|u\|^4_\frac{12}{5};
\end{equation}
\item if $u$ is a radial function then $\phi_u$ is radial, too. \label{phi6}
\end{enumerate}
\end{lemma}
%
%
%

Now we use the previous lemma to deduce assumptions
\ref{r1}-\ref{r6}.\\
Hypothesis \ref{r1} is obvious.\\
Since
\[
R'(u)[u]=\irt \phi_u u^2,
\]
(see for example \cite{BF}), then \ref{rd'} is again a consequence of \eqref{eq:phiq}.
\\
We pass to check \ref{r4}. Suppose that
    \begin{equation*}
        u_j \rightharpoonup u \hbox{ weakly in }\HTr.
    \end{equation*}
By compact embedding we deduce that
    \begin{equation*}
        u_j\to u \hbox{ in } L^{\frac {12} 5}(\RT)
    \end{equation*}
and then, by continuity,
    \begin{equation*}
        \phi_{u_j}\to \phi_u \hbox{ in } \D.
    \end{equation*}
Since $R'(u)[v]=\irn \phi_u u v$, we
have that
    \begin{align*}
        \limsup_j R'(u_j)[u-u_j]&=\limsup_j \irt \phi_{u_j}u_j(u-u_j)
        \\&\le C \limsup_j  \|\phi_{u_j}\|_{\D}\|u_j\|_{\frac{12}5} \|u-u_j\|_{\frac
        {12}5}=0.
    \end{align*}
Now in order to verify \ref{r5}, we consider $u\in\HT,$ $u\neq 0$
and the rescaled function $u_t$. We compute
    $$R(u_{t})=\frac 1 4 \irt \phi_{u_t}u_t^2=\frac {t^5} 4\irt \phi_u u^2=t^5R(u)$$
so \ref{r5} holds true for $\a =5$.\\
Finally \ref{r6} follows from
\ref{phi6} of Lemma \ref{le:prop}.

\subsection{The elliptic Kirchhoff equation}

In this subsection we treat the semilinear
perturbation of the Kirchhoff equation
    \begin{equation}\tag{$\mathcal K$}\label{eq:smkir}
        - \left(p+q\irn |\n u|^2\right) \Delta u=g(u)  \qquad \hbox{ in }\RN,
    \end{equation}
where $p>0$ and $g$ satisfies \ref{g1}-\ref{g4}. Arguing as in \cite{A,BL1}, without loss of generality, we can suppose that $g$ satisfies \ref{g3'}. We find the solution to \eqref{eq:smkir} as the
critical points of the functional
    \begin{equation*}
        I_q(u)= \frac 1 2 \left( p +\frac q 2\irn|\n u|^2\right)\irn|\n u|^2-\irn
        G(u).
    \end{equation*}
It is easy to see that $I_q$ is of the type \eqref{eq:funq}, where
$R(u)=\frac 1 4 \left(\irn |\n u|^2\right)^2.$\\
Assumptions \ref{r1}-\ref{rd'} are trivially satisfied as we can see
by straight computations.\\
As to \ref{r4}, suppose that
$u_j\rightharpoonup u$ weakly in $\Hr.$ By weak lower semicontinuity, we
know that
    \begin{equation*}
        \irn |\n u|^2\le \liminf_j \irn |\n u_j|^2,
    \end{equation*}
and then
    \begin{align*}
        \limsup_j R'(u_j)[u-u_j]&=\limsup_j\irn  |\n u_j |^2 \cdot\irn \n u_j\cdot \n
        (u-u_j)\\
        &\le \limsup_j \irn |\n u_j|^2\cdot\limsup_j \irn  \n u_j \cdot \n
        (u-u_j)\\
        &\le \limsup_j \irn |\n u_j|^2\cdot\left(\limsup_j\irn\n u_j\cdot\n
        u-\liminf_j\irn |\n u_j|^2 \right)\\
        &\le \limsup_j \irn |\n u_j|^2\cdot\left(\irn|\n u|^2-\liminf_j\irn |\n u_j|^2
        \right)\le 0.
    \end{align*}
By a simple computation, we have that
    $$R(u_t)=\frac 1 4 \left(\irn |\n u_t|^2\right)^2= \frac
    {t^{2(N-2)}}{4}\left(\irn |\n u|^2\right)^2=t^{2(N-2)} R(u),$$
and then also \ref{r5} is satisfied.\\
Finally by a simple change of
variable it can be proved that for any $g\in O(N)$ we have
    $$R(u(g x))=\frac 1 4 \left(\irn |\n u(g x)|^2\right)^2=\frac 1 4 \left(\irn |\n u(x)|^2\right)^2=R(u).$$

\begin{remark}
Let us observe that we can easily apply Theorem \ref{main} also to a sort of linear combination of the Schr\"odinger-Maxwell equation with the Kirchhoff one, namely we can find multiple critical points of the functional
\begin{equation*}
I_q(u)= \frac 12 \irt |\n u|^2
+\frac q4 \left[\l_1 \irt \left(\frac 1 {|x|}\ast u^2\right)u^2
+\l_2\left( \irt |\n u|^2 \right)^2\right]
-\irt G(u),
\end{equation*}
with $\l_1,\l_2\in \R_+$ and $q$ sufficiently small.
\end{remark}

%
%
%
%

%
%
%

\end{document}